\documentclass[a4paper, fontsize=12pt]{article}

\renewenvironment{abstract}
               {\list{}{\rightmargin\leftmargin}%
                \item[\hspace*{1cm}\small\textbf{Abstract ---}]\relax}
               {\endlist}

\usepackage[top=0.8in, bottom=0.8in, left=0.8in, right=0.8in]{geometry}
\usepackage[fleqn]{amsmath}
\usepackage{amssymb}
\usepackage{enumerate}
\usepackage{amsthm,bm}
\usepackage{appendix}
\usepackage[russian,english]{babel}
\usepackage{subfigure}
\usepackage{float}
\usepackage[pdftex]{graphicx}
\usepackage{framed}
\usepackage{hyperref}
\usepackage{sidecap}
\usepackage{mathrsfs}

\sidecaptionvpos{figure}{t}
\newtheorem{Theorem}{Theorem}

\newtheorem{Corollary}[Theorem]{Corollary}

\newtheorem{Definition}[Theorem]{Definition}

\newtheorem{Proposition}[Theorem]{Proposition}

\begin{document}

\title{\textbf{Relative consistency of Set Matrix Theory with ZF}}

\author{Marcoen J.T.F. Cabbolet\\
        \small{\textit{Center for Logic and Philosophy of Science, Vrije Universiteit Brussel}}\footnote{email: Marcoen.Cabbolet@vub.be}
        }
\date{}

\maketitle
\vfill
\begin{abstract}
Set Matrix Theory (SMT) has been introduced in \emph{Log. Anal.} \textbf{225}: 59-82 (2014) as a generalization of ZF, in which matrices constructed from sets are treated as urelements, that is, as objects that are not sets but that can be elements of sets. Here we prove that SMT is relatively consistent with ZF.
\end{abstract}

\noindent In the context of axiomatic set theories, we may define \emph{urelements} as objects that are not sets but that can be elements of sets \cite{Mendelson}. ZF is then an axiomatic set theory without urelements: every object is a set. Set Matrix Theory (SMT), on the other hand, has been introduced in \cite{SMT} as a generalization of ZF in which matrices constructed in a finite number of steps from a finite number of sets are treated as \emph{urelements}. These urelements, however, are not ``primordial'' in the sense of ``not constructible'': in the framework of SMT, the urelement (set matrices) are constructed from sets. That being said, Shoenfield defined urelements as objects that are not sets \emph{and} that do not involve sets in their construction \cite{Shoenfield}: from that perspective, the set matrices in the universe of SMT cannot be called urelements. Nevertheless, in this paper we will use the term `urelements' for the set matrices in the framework of SMT: those who object against this usage may replace it by `objects sui generis'.

The language $L_{SMT}$ of SMT is the language $L_{ZF}$ of ZF, with the constant $\emptyset$ and with (labeled) Roman variables $x, x_1, \ldots, y, y_1, \ldots$ varying over sets, extended with
\begin{enumerate}[(i)]
  \item (labeled) Greek variables $\alpha,\alpha_1, \ldots, \beta, \beta_1, \ldots$ varying over set matrices;
  \item for each pair of positive integers $m$ and $n$, an $mn$-ary function symbols $f_{m\times n}$, with matrix notation
  $\left [ \begin{array}{ccc} t_{11} & \cdots & t_{1n} \\ \vdots &   & \vdots \\ t_{m1} & \cdots & t_{mn} \end{array} \right ]$ for $f_{m\times n}(t_{11}, \ldots, t_{mn})$, where the $t_{ij}$'s are terms of $L_{SMT}$.
\end{enumerate}
A constructive \textbf{set matrix axiom schema} states that given any $m\cdot n$ set matrices $\alpha_{11}, \ldots, \alpha_{mn}$, there exists an $m\times n$ set matrix $\beta$ whose entries are the $\alpha_{ij}$'s:
\begin{equation}\label{eq:SetMatrixAxiom}
\forall \alpha_{11} \cdots \forall \alpha_{mn} \exists \beta: \beta = \left [ \begin{array}{ccc} \alpha_{11} & \cdots & \alpha_{1n} \\ \vdots &   & \vdots \\ \alpha_{m1} & \cdots & \alpha_{mn} \end{array} \right ]
\end{equation}
The \textbf{reduction axiom} of SMT identifies a $1\times 1$ set matrix with its entry:
\begin{equation}\label{eq:ReductionAxiom}
  \forall x : [x] = x
\end{equation}
This way, quantification over all set matrices implies quantification over all sets. Likewise, the \textbf{omission axiom schema} identifies a $1\times 1$ set matrix whose entry is an $m\times n$ set matrix with its entry:
\begin{equation}\label{eq:OmissionAxiom}
\forall \alpha_{11} \cdots \forall \alpha_{mn} : \left [ \left [ \begin{array}{ccc} \alpha_{11} & \cdots & \alpha_{1n} \\ \vdots &   & \vdots \\ \alpha_{m1} & \cdots & \alpha_{mn} \end{array} \right ] \right ]= \left [ \begin{array}{ccc} \alpha_{11} & \cdots & \alpha_{1n} \\ \vdots &   & \vdots \\ \alpha_{m1} & \cdots & \alpha_{mn} \end{array} \right ]
\end{equation}
\newpage
The \textbf{division axiom schema} of SMT then states that set matrices of different sizes are not identical. This schema consists of a subschema \eqref{eq:DivisionAxiomSets} that introduces equalities between sets and various urelements, and a subschema \eqref{eq:DivisionAxiomMatrices} that introduces inequalities between various urlements:
\begin{subequations}
\begin{equation}\label{eq:DivisionAxiomSets}
\forall x \forall \alpha_{11} \cdots \forall \alpha_{mn}: x \neq \left [ \begin{array}{ccc} \alpha_{11} & \cdots & \alpha_{1n} \\ \vdots &   & \vdots \\ \alpha_{m1} & \cdots & \alpha_{mn} \end{array} \right ] \ \ \ (m\cdot n \geq 2)
\end{equation}
\begin{equation}\label{eq:DivisionAxiomMatrices}
\forall \alpha_{11} \cdots \forall \alpha_{mn}\forall \beta_{11} \cdots \forall \beta_{pq}: \left [ \begin{array}{ccc} \alpha_{11} & \cdots & \alpha_{1n} \\ \vdots &   & \vdots \\ \alpha_{m1} & \cdots & \alpha_{mn} \end{array} \right ] \neq  \left [ \begin{array}{ccc} \beta_{11} & \cdots & \beta_{1q} \\ \vdots &   & \vdots \\ \beta_{p1} & \cdots & \beta_{pq} \end{array} \right ] \ \ \ (1,1) \neq (m,n) \neq (p,q) \neq (1,1)
\end{equation}
\end{subequations}
And the \textbf{epsilon axiom schema} states that an $m\times n$ set matrix with $m\cdot n \geq 2$ has no elements in the sense of the $\in$-relation:
\begin{equation}\label{eq:EpsilonAxiom}
\forall \alpha_{11} \cdots \forall \alpha_{mn} \forall \beta: \beta \not\in \left [ \begin{array}{ccc} \alpha_{11} & \cdots & \alpha_{1n} \\ \vdots &   & \vdots \\ \alpha_{m1} & \cdots & \alpha_{mn} \end{array} \right ] \ \ \ (m\cdot n \geq 2)
\end{equation}
It is therefore that such $m\times n$ set matrices can be viewed as \emph{urelements}. The \textbf{extensionality axiom schema for set matrices} then states that two $m\times n$ set matrices are equal if their entries are equal:
\begin{equation}\label{eq:SetMatricesExt}
\forall \alpha_{11} \cdots \forall \alpha_{mn} \forall \beta_{11} \cdots \forall \beta_{mn}  : \left [ \begin{array}{ccc} \alpha_{11} & \cdots & \alpha_{1n} \\ \vdots &   & \vdots \\ \alpha_{m1} & \cdots & \alpha_{mn} \end{array} \right ] = \left [ \begin{array}{ccc} \beta_{11} & \cdots & \beta_{1n} \\ \vdots &   & \vdots \\ \beta_{m1} & \cdots & \beta_{mn} \end{array} \right ] \Leftrightarrow \alpha_{11} = \beta_{11} \wedge \ldots \wedge \alpha_{mn} = \beta_{mn}
\end{equation}
As to the set-theoretical axioms of SMT, these are generalizations of the axioms of ZF, since set matrices that are not sets can still be elements of sets. E.g. the \textbf{generalized extensionality axiom for sets} becomes
\begin{equation}\label{eq:GeneralizedEXT}
  \forall x \forall y : x = y \Leftrightarrow \forall \alpha (\alpha \in x \Leftrightarrow \alpha\in y)
\end{equation}
The generalizations of the constructive set-theoretical axioms of ZF are then the following:
\begin{gather}
\exists x : x = \emptyset \wedge \forall \alpha (\alpha \not \in x)\label{eq:GeneralizedEMPTY} \\
\forall x \exists y \forall \alpha : \alpha \in y \Leftrightarrow \alpha \in x \wedge \Phi(\alpha)  \\
\forall \alpha \forall \beta \exists x \forall \gamma : \gamma \in x \Leftrightarrow \gamma = \alpha \vee \gamma = \beta \\
\forall x : \forall \alpha (\alpha \in x \Rightarrow \exists u (u = \alpha )) \Rightarrow \exists y \forall \beta (\beta \in y \Leftrightarrow \exists z (z \in x \wedge \beta \in z)) \\
\forall x \exists y \forall \alpha : \alpha \in y \Leftrightarrow \exists u (u = \alpha \wedge u \subseteq x) \\
\exists x : \emptyset \in x \wedge \forall y (y \in x \Rightarrow \{y \} \in x) \\
\forall x : \forall \alpha(\alpha \in x \Rightarrow \exists!\beta : \Phi(\alpha, \beta)) \Rightarrow \exists y \forall \zeta (\zeta\in y \Leftrightarrow \exists \gamma : \gamma \in x \wedge \Phi(\gamma, \zeta))
\end{gather}
In addition, SMT has the \textbf{set of set matrices axiom schema}:
\begin{equation}\label{eq:SetOfSetMatrices}
  \forall x \exists y \forall \alpha : \alpha \in y \Leftrightarrow
  \exists \beta_{11} \ldots \exists \beta_{mn}(\alpha = \left [ \begin{array}{ccc} \beta_{11} & \cdots & \beta_{1n} \\ \vdots &   & \vdots \\ \beta_{m1} & \cdots & \beta_{mn} \end{array} \right ] \wedge
\beta_{11} \in x \wedge \ldots \wedge \beta_{mn} \in x)
\end{equation}
These are all the constructive axioms of SMT. The axiom of regularity of SMT plays no role in the present discussion.\\
\ \\
We now want to prove that SMT is relatively consistent with ZF, that is,
\begin{equation}
  \left[\ \not\vdash_{ZF}\bot\ \right] \longrightarrow \left[\ \not\vdash_{SMT}\bot\ \right]
\end{equation}
For that matter, we first construct the theory ZFM by extending $L_{ZF}$ with function symbols $f_{m\times n}$ and adding an axiom schema $\Psi(f_{m\times n})$ to ZF, so that for positive integers $m$ and $n$ the axiom $\Psi(f_{m\times n})$ defines $m\times n$ set matrices in the framework of ZF (Def. \ref{def:SetMatricesZF}). Thereafter we construct the theory SMT$^-$ by removing the division axiom schema \eqref{eq:DivisionAxiomSets}and the epsilon axiom schema \eqref{eq:EpsilonAxiom} from SMT, and we show that there is an interpretation of SMT$^-$ in ZFM (Prop. \ref{prop:interpretation}). We then prove that SMT is relatively consistent with SMT$^-$ (Prop. \ref{prop:SMTrcwSMT-}). From these results we derive that SMT is relatively consistent with ZF (Prop. \ref{prop:SMTrcwZF}).
\newpage
\begin{Definition}[Set matrices in ZF]\label{def:SetMatricesZF}\ \\
In the framework of ZF, we add function symbols $f_{m\times n}$ to the language $L_{ZF}$; a \textbf{set matrix} $\left( \begin{array}{ccc} x_{11} & \cdots & x_{1n} \\ \vdots &   & \vdots \\ x_{m1} & \cdots & x_{mn} \end{array} \right)$ is then a notation for an object $f_{m\times n}(x_{11},\ldots, x_{mn})$, given by the defining axiom for $f_{m\times n}$:
\begin{subequations}\label{eq:MatrixModel}
\begin{gather}\label{eq:1x1SetMatricesZFM}
\forall x_{11} : (x_{11}) = x_{11} \\
\forall x_{11} \ldots  \forall x_{mn} : \left( \begin{array}{ccc} x_{11} & \cdots & x_{1n} \\ \vdots &   & \vdots \\ x_{m1} & \cdots & x_{mn} \end{array} \right)= \{((1,1), x_{11}), \ldots, ((m,n), x_{m,n})\} \ \ \ (m\cdot n> 1)\label{eq:mxnSetMatricesZFM}
\end{gather}
\end{subequations}
We now also add (labeled) Greek variables $\alpha,\alpha_1, \ldots, \beta, \beta_1, \ldots$, ranging over set matrices defined by the schema \eqref{eq:MatrixModel}, to the language. ZFM is then the theory obtained by adding the defining axioms \eqref{eq:MatrixModel} as an axiom schema to ZF. (Note that the Greek variables are part of the language $L_{ZFM}$ of ZFM, but do not occur in the axioms of ZFM.)\hfill$\blacksquare$
\end{Definition}

\begin{Proposition}\label{prop:ZFMrcwZF}
ZFM is relatively consistent with ZF.\hfill$\blacksquare$
\end{Proposition}
\paragraph{Proof:} \ \\
ZFM is a \emph{definitional extension} of ZF, that is, an extension of ZF with the defining axioms \eqref{eq:MatrixModel} for the function symbols $f_{m\times n}$. But then ZFM is a conservative extension of ZF, cf. \cite{Shoenfield}. That means that for every nonlogical axiom $\Phi$ of ZF and for every sentence $\Psi$ in $L_{ZF}$ we have
\begin{subequations}
\begin{gather}
\vdash_{ZFM}\Phi\\
\left[\ \vdash_{ZFM}\Psi\ \right] \longrightarrow \left[\ \vdash_{ZF}\Psi\ \right]
\end{gather}
\end{subequations}
Being a conservative extension of ZF, ZFM is relatively consistent with ZF---cf. \cite{Shoenfield}. \hfill$\blacksquare$

\begin{Proposition}
For any sentence $\forall x\Psi(x)$ with an occurrence of the variable $x$ in $L_{ZFM}$ we have
\begin{equation}
\vdash_{ZFM}\ \forall x\Psi(x) \Leftrightarrow [\alpha/ x]\forall x\Psi(x)
\end{equation}
where $[\alpha/ x]\forall x\Psi(x)$ is the sentence obtained by replacing every occurrence of $x$ in $\forall x\Psi(x)$ by $\alpha$.\hfill$\blacksquare$
\end{Proposition}
\paragraph{Proof:} \ \\
Since by axiom schema \eqref{eq:MatrixModel} set matrices are just sets in the framework of ZFM, quantification over all sets in the framework of ZFM implies quantification over all set matrices. And vice versa: since by axiom \eqref{eq:1x1SetMatricesZFM} a $1\times 1$ matrix with a set-valued entry is identical to that set-valued entry, quantification over all set matrices---this implies quantification over all $1\times 1$ set matrices---implies quantification over all sets.\hfill$\blacksquare$

\begin{Proposition}
For every well-formed formula $\Psi$ in $L_{SMT}$ there is an identical well-formed formula $\Psi^\prime$ in $L_{ZFM}$.\hfill$\blacksquare$
\end{Proposition}
\paragraph{Proof:} \ \\
$L_{SMT}$ is $L_{ZF}$ extended with (labeled) Greek variables $\alpha,\alpha_1, \ldots, \beta, \beta_1, \ldots$ ranging over set matrices, and with function symbols $f_{m\times n}$. By Def. \ref{def:SetMatricesZF}, $L_{ZFM}$ is $L_{ZF}$ extended with those same Greek variables and function symbols. Ergo, $L_{SMT}$ and $L_{ZFM}$ are the same formal languages, which means that for every well-formed formula $\Psi$ in $L_{SMT}$ there is an identical well-formed formula $\Psi^\prime$ in $L_{ZFM}$.\hfill$\blacksquare$\\
\ \\
Of course, in ZFM we have the mathematical axioms \eqref{eq:mxnSetMatricesZFM} by which $f_{m\times n}(x_{11},\ldots, x_{mn})$ is identical to $\{((1,1), x_{11}), \ldots, ((m,n), x_{m,n})\}$. We don't have these axioms in SMT. But that doesn't matter: what matters is that we have identical well-formed formulas in $L_{SMT}$ and $L_{ZFM}$. With that, we have a very simple translation from $L_{SMT}$ to $L_{ZFM}$, by which a formula in $L_{SMT}$ with occurrences of matrices with square brackets `[' and `]' is translated to the formula in $L_{ZFM}$ with occurrences of matrices with round brackets `(' and `)' instead. For example, $[[\emptyset\ \ \emptyset]\ \emptyset\ ]$ is a notation for the term $f_{1\times2}(f_{1\times2}(\emptyset, \emptyset), \emptyset)$ of $L_{SMT}$ and $((\emptyset\ \ \emptyset)\ \emptyset\ )$ is a notation for the (identical) term $f_{1\times2}(f_{1\times2}(\emptyset, \emptyset), \emptyset)$ of $L_{ZFM}$; so, a well-formed formula $\Psi([[\emptyset\ \ \emptyset]\ \emptyset\ ])$ in $L_{SMT}$ with an occurrence of $[[\emptyset\ \ \emptyset]\ \emptyset\ ]$ is translated to the well-formed formula $\Psi(((\emptyset\ \ \emptyset)\ \emptyset\ ))$ in $L_{ZFM}$ obtained by replacing every occurrence of $[[\emptyset\ \ \emptyset]\ \emptyset\ ]$ in $\Psi([[\emptyset\ \ \emptyset]\ \emptyset\ ])$ by $((\emptyset\ \ \emptyset)\ \emptyset\ )$.

\begin{Proposition} Let $\Psi$ be an axiom of SMT, which is not an axiom of the epsilon axiom schema \eqref{eq:EpsilonAxiom} nor an axiom of the division axiom scheme \eqref{eq:DivisionAxiomSets}, and let $\Psi^\prime$ be the corresponding well-formed formula in $L_{ZFM}$. Then
\begin{equation}
\left[\ \vdash_{SMT} \Psi\ \right] \longrightarrow \left[\ \vdash_{ZFM} \Psi^\prime\ \right]
\end{equation}
$\blacksquare$
\end{Proposition}

\paragraph{Proof:}\ \\
\noindent Let $m,n$ be finite Von Neumann ordinals in the framework of ZF, with $mn>1$; then it is a theorem of ZF that for any $m\cdot n$ sets $x_{ij}$ there is precisely one function $f$ on the cartesian product $\{1, \ldots, m\}\times\{1, \ldots, n\}$ such that $f$ maps $(i,j)$ to $x_{ij}$:
\begin{equation}
\forall x_{11} \ldots  \forall x_{mn} \exists! f : f = \{((1,1), x_{11}), \ldots, ((m,n), x_{mn})\}
\end{equation}
So, as a corollary we have in ZFM that
\begin{equation}\label{eq:MatrixCol}
\forall \alpha_{11} \cdots \forall \alpha_{mn} \exists! y: y = \left( \begin{array}{ccc} \alpha_{11} & \cdots & \alpha_{1n} \\ \vdots &   & \vdots \\ \alpha_{m1} & \cdots & \alpha_{mn} \end{array} \right)\ \ \ (m\cdot n>2)
\end{equation}
Together with Eq. \eqref{eq:1x1SetMatricesZFM}, this proves that each axiom of the set matrix axiom schema \eqref{eq:SetMatrixAxiom}, identically interpreted in $L_{ZFM}$, is a theorem of ZFM. Furthermore, the interpretations of the reduction axiom of SMT, Eq. \eqref{eq:ReductionAxiom}, and the omission axiom schema, Eqs. \eqref{eq:OmissionAxiom}, in $L_{ZFM}$ follow directly from the axiom schema \eqref{eq:MatrixModel} in the framework of ZFM.

Proceeding, the interpretations of the division axiom schema for set matrices of SMT, Eqs. \eqref{eq:DivisionAxiomMatrices}, and the extensionality axiom schema for set matrices of SMT, Eqs. \eqref{eq:SetMatricesExt}, can be derived in the framework of ZFM directly from the axiom schema \eqref{eq:MatrixModel} and the extensionality axiom for sets of ZF. This covers all the matrix-theoretical axioms of SMT: we have herewith proven that for every matrix-theoretical axiom of SMT, which is not an epsilon axiom \eqref{eq:EpsilonAxiom} nor a division axiom \eqref{eq:DivisionAxiomSets}, there is an identical well-formed formula in $L_{ZFM}$ that is a theorem in the framework of ZFM.

As to the set-theoretical axioms of SMT, Eqs. \eqref{eq:GeneralizedEXT}-\eqref{eq:SetOfSetMatrices}, the crux is that these are all simple theorems of ZF---the proof is left as an exercise for the reader. That proves the proposition: every axiom $\Psi$ of SMT, which is not a division axiom \eqref{eq:DivisionAxiomSets} nor an epsilon axiom \eqref{eq:EpsilonAxiom}, is, identically interpreted in $L_{ZFM}$, a theorem $\Psi^\prime$ of ZFM.\hfill$\blacksquare$

\begin{Proposition}\label{prop:interpretation}
Let SMT$^{-}$ be the subtheory of SMT obtained by removing the epsilon axiom schema \eqref{eq:EpsilonAxiom} and the division axioms of the scheme \eqref{eq:DivisionAxiomSets}; then there is an interpretation of SMT$^{-}$ in ZFM.\hfill$\blacksquare$
\end{Proposition}

\paragraph{Proof:}\ \\
We have proven that every nonlogical axiom of SMT$^{-}$, identically interpreted in $L_{ZFM}$, is a theorem of ZFM. \hfill$\blacksquare$

\begin{Corollary}
SMT$^{-}$ is relatively consistent with ZFM.\hfill$\blacksquare$
\end{Corollary}

\begin{Proposition}\label{prop:SMTrcwSMT-}
SMT is relatively consistent with SMT$^{-}$.\hfill$\blacksquare$
\end{Proposition}

\paragraph{Proof:}\ \\
Observe that SMT is just the theory SMT$^{-}\ \cup \{\Psi(f_{1\times 2}), \Psi(f_{2\times 1}), \Psi(f_{1\times 3}), \ldots, \Phi(f_{1\times 2}), \Phi(f_{2\times 1}), \Phi(f_{1\times 3}), \ldots, \}$ obtained by extending SMT$^{-}$ with the epsilon axioms $\Psi(f_{m\times n})$ of schema \eqref{eq:EpsilonAxiom} and with the division axioms $\Phi(f_{m\times n})$ of schema \eqref{eq:DivisionAxiomSets}. Assuming that SMT$^{-}$ is consistent, to prove Prop. \ref{prop:SMTrcwSMT-} we need to prove that we cannot derive the negations of the epsilon axioms $\Psi(f_{m\times n})$ from the axioms of SMT$^-$, and that we cannot derive the negations of the division axioms $\Phi(f_{m\times n})$ from SMT$^{-} \cup \{\Psi(f_{1\times 2}), \Psi(f_{2\times 1}), \Psi(f_{1\times 3}), \ldots, \}$. In other words, to prove Prop. \ref{prop:SMTrcwSMT-}, we will prove
\begin{subequations}
\begin{gather}\label{eq:NotNotEpsilonAxioms}
  \not\vdash_{SMT^-}\neg \Psi(f_{m\times n})\\
  \not\vdash_{SMT^- \cup \{\Psi(f_{1\times 2}), \Psi(f_{2\times 1}), \Psi(f_{1\times 3}), \ldots, \}}\neg \Phi(f_{m\times n})\label{eq:NotNotDivisionAxioms}
\end{gather}
\end{subequations}
So, note that for the proof of Prop. \ref{prop:SMTrcwSMT-} we stay in the framework of SMT$^-$: we do not use the translation from $L_{SMT}$ to $L_{ZFM}$.

That being said, let us say that the position of the term $t_1$ to the left of the $\in$-symbol in an atomic formula $t_1\in t_2$ is the \textbf{element position}, and that the position of the term $t_2$ to the right of the $\in$-symbol in an atomic formula $t_1\in t_2$ is the \textbf{set position}. We then observe that there is no axiom of SMT$^{-}$ which contains a subformula of the type $t_1\in t_2$ with an occurrence of a (Greek) variable ranging over set matrices on the set position of the $\in$-relation. But that means that for every epsilon axiom $\Psi(f_{m\times n})$ of the schema \eqref{eq:EpsilonAxiom}, we cannot derive $\neg\Psi(f_{m\times n})$ from SMT$^{-}$:
\begin{equation}
  \not\vdash_{SMT^-}\neg\Psi(f_{m\times n})
\end{equation}
But that is precisely \eqref{eq:NotNotEpsilonAxioms}.
\newpage
\noindent Proceeding with the same method, we observe
\begin{enumerate}[(i)]
  \item that there is no axiom of the theory SMT$^{-} \cup \{\Psi(f_{1\times 2}), \Psi(f_{2\times 1}), \Psi(f_{1\times 3}), \ldots, \}$, obtained by extending SMT$^{-}$ with the epsilon axiom schema $\Psi(f_{1\times 2}), \Psi(f_{2\times 1}), \Psi(f_{1\times 3}), \ldots$, which contains an identity between a set and a set matrix constructed from more than one set;
  \item that we cannot derive an identity between a set and a set matrix constructed from more than one set from the extensionality axiom for sets \eqref{eq:GeneralizedEXT} or from the extensionality axioms for set matrices \eqref{eq:SetMatricesExt}, since these can only yield identities between sets or identities between set matrices constructed from the same number of sets.
\end{enumerate}
But that means that for every division axiom $\Phi(f_{m\times n})$ of the schema \eqref{eq:EpsilonAxiom}, we cannot derive $\neg\Phi(f_{m\times n})$ from SMT$^{-} \cup \{\Psi(f_{1\times 2}), \Psi(f_{2\times 1}), \Psi(f_{1\times 3}), \ldots, \}$:
\begin{equation}
  \not\vdash_{SMT^- \cup \{\Psi(f_{1\times 2}), \Psi(f_{2\times 1}), \Psi(f_{1\times 3}), \ldots, \}}\neg\Phi(f_{m\times n})
\end{equation}
But that is precisely \eqref{eq:NotNotDivisionAxioms}. Ergo, SMT is relatively consistent with SMT$^{-}$. \hfill$\blacksquare$

\begin{Proposition}\label{prop:SMTrcwZF}
SMT is relatively consistent with ZF.\hfill$\blacksquare$
\end{Proposition}

\paragraph{Proof:}\ \\
We have proven
\begin{enumerate}[(i)]
  \item that SMT is relatively consistent with SMT$^{-}$;
  \item that SMT$^{-}$ is relatively consistent with ZFM;
  \item that ZFM is relatively consistent with ZF.
\end{enumerate}
Ergo, by transitivity of relative consistency, SMT is relatively consistent with ZF.\hfill$\blacksquare$\\
\ \\
This shows that SMT can be used as a foundational theory whenever it is desirable to treat set matrices as \emph{objects sui generis}. One just has to be careful with the definition of a transitive set, which is used for the definition of ordinals. There are two common definitions of a transitive set:
\begin{enumerate}[(i)]
  \item a set $x$ is \textbf{transitive} if $\forall\alpha\forall\beta: \alpha\in x \wedge \beta\in \alpha\Rightarrow \beta\in x$ \cite{Shoenfield,Bernays}
  \item a set $x$ is \textbf{transitive} if every element $y$ of $x$ is a subset of $x$ \cite{Shoenfield}
\end{enumerate}
In \cite{Shoenfield} it is stated that these are equivalent, but in the framework of SMT that is not the case. For example, the set $x = \{[\emptyset\ \ \emptyset], \{[\emptyset\ \ \emptyset]\}\}$ satisfies definition (i) because every element $\alpha$ of every set $y$ in $x$ is in $x$. But it does not satisfy definition (ii), because the object $[\emptyset\ \ \emptyset]$ is not a set hence not a subset of $x$. Now what we want in axiomatic set theory is that all elements of a transitive set are sets, so we can avoid weird transitive sets like the set $x = \{[\emptyset\ \ \emptyset], \{[\emptyset\ \ \emptyset]\}\}$ by either using definition (ii) or by replacing definition (i) by the following:
\begin{enumerate}[(i)]
  \item[(iii)] a set $x$ is \textbf{transitive} if every element $\alpha\in x$ is either the empty set or there is a $\beta$ such that $\beta\in\alpha$, and if every element $\gamma$ of every $\alpha\in x$ is in $x$.
\end{enumerate}
In the framework of SMT, definitions (ii) and (iii) are equivalent. One may then proceed by defining ordinals as transitive sets well-ordered by $\in$.

\end{document}